\documentclass[conference]{IEEEtran}
\IEEEoverridecommandlockouts
\usepackage{cite}
\usepackage{amsmath,amssymb,amsfonts}
\usepackage{algorithmic}
\usepackage{graphicx}
\usepackage{textcomp}
\usepackage{xcolor}
\def\BibTeX{{\rm B\kern-.05em{\sc i\kern-.025em b}\kern-.08em
    T\kern-.1667em\lower.7ex\hbox{E}\kern-.125emX}}

\usepackage[ruled]{algorithm2e}

\newtheorem{theorem}{Theorem}

\newcommand{\calC}{\mathcal{C}}
\newcommand{\calG}{\mathcal{G}}
\newcommand{\calN}{\mathcal{N}}
\newcommand{\calK}{\mathcal{K}}
\newcommand{\calL}{\mathcal{L}}
\newcommand{\calF}{\mathcal{F}}
\newcommand{\calW}{\mathcal{W}}
\newcommand{\calD}{\mathcal{D}}
\newcommand{\calB}{\mathcal{B}}
\newcommand{\calR}{\mathcal{R}}

\newcommand{\ibf}{\mathbf{i}}
\newcommand{\setC}{\mathbb{C}}

\DeclareMathOperator{\conv}{\text{conv}}

\begin{document}

\title{On the Tightness of the Lagrangian Dual Bound for Alternating Current Optimal Power Flow\thanks{This material is based upon work supported by the U.S. Department of Energy, Office of Science, Advanced Scientific Computing Research, under Contract DE-AC02-06CH11357.  We also acknowledge support from the U.S. NSF under award 1832208. We gratefully acknowledge the computing resources provided on Bebop, a high-performance computing cluster operated by the Laboratory Computing Resource Center at Argonne National Laboratory.}
}

\author{\IEEEauthorblockN{Weiqi Zhang}
\IEEEauthorblockA{
\textit{University of Wisconsin-Madison}\\
Madison, WI, USA \\
wzhang483@wisc.edu}
\and
\IEEEauthorblockN{Kibaek Kim}
\IEEEauthorblockA{
\textit{Argonne National Laboratory}\\
Lemont, IL, USA \\
kimk@anl.gov}
\and
\IEEEauthorblockN{Victor M. Zavala}
\IEEEauthorblockA{
\textit{University of Wisconsin-Madison}\\
Madison, WI, USA \\
zavalatejeda@wisc.edu}
}

\maketitle

\begin{abstract}
We study tightness properties of a Lagrangian dual (LD) bound for the nonconvex alternating current optimal power flow (ACOPF) problem. We show an LD bound that can be computed in a parallel, decentralized manner. Specifically, the proposed approach partitions the network into a set of subnetworks, dualizes the coupling constraints (giving the LD function), and maximizes the LD function with respect to the dual variables of the coupling constraints (giving the desired LD bound). The dual variables that maximize the LD are obtained by using a proximal bundle method. We show that the bound is less tight than the popular semidefinite programming relaxation but as tight as the second-order cone programming relaxation. We demonstrate our developments using PGLib-OPF test instances. 
\end{abstract}

\begin{IEEEkeywords}
alternating current optimal power flow, Lagrangian duality, distributed optimization
\end{IEEEkeywords}

\section{Introduction}
%
%
%
%
We study the solution of alternating current optimal power flow (ACOPF) problem. Finding a globally optimal solution for this nonconvex problem is an NP-hard problem. Local solution algorithms have been widely used and accepted for finding a local optimal solutions. However, the quality of such local solutions cannot be assessed without having a lower bound. Several relaxation schemes have been proposed to obtain a lower bound of the solution. Popular relaxation schemes include the second-order cone programming (SOCP) relaxation \cite{soc_form}, the semidefinite programming (SDP) relaxation (e.g., \cite{sdp_original}), and QC relaxation~\cite{coffrin2015qc}. 

Here, we propose a new relaxation approach and study the tightness properties of its associated lower bound. The approach partitions the network into several subnetworks and applies a Lagrangian relaxation of the coupling constraints. This relaxation/dualization approach is used to formulate the so-called Lagrangian dual (LD) function and the lower bound is bound by maximizing the LD function with respect to the the dual variables of the coupling constraints. This decomposition approach, also commonly known as LD decomposition, has been successfully applied to tackle large-scale nonconvex optimization problems arising in different application domains. The success of this approach can be attributed to the observation that LD bounds are in general quite tight  (e.g.,~\cite{kim2018algorithmic,kim2018temporal}). Here, we prove that the LD bound provided by this approach is tighter than that of SOCP but more relaxed than that of SDP; however, a key advantage of our approach is that the LD problem can be solved in a parallel/decentralized manner (thus enabling handling of large-scale networks). 
Note that our relaxation approach can be generalized to other variants of ACOPF (e.g. multi-period and/or stochastic ACOPF) by exploiting the problem structures.


From an operations perspective, an important observation is that the proposed approach can be interpreted as a decentralized architecture that aims to find an approximate solution of the ACOPF problem. In particular, future electric grids face many challenges, including increasing penetrations of distributed energy resources, cybersecurity threats, and data privacy (e.g.,~\cite{ryu2021privacy}), which necessitate the transition from centralized to decentralized methods \cite{dopf_survey}. 
However, most existing decentralized algorithms such as the alternating direction method of multipliers (e.g.,~\cite{ admm_full_sun,kim2021leveraging}) and the analytical target cascading (e.g.,~\cite{dormohammadi2013exponential}) do not provide a lower bound for the problem (and thus optimality cannot be easily assessed). 
In this work, we derive a proximal bundle method (PBM) to solve the LD problem and provide convergence guarantees for such scheme; importantly, the lower bound obtained with this approach can be used for assessing the quality of the approximate solution. 


The contribution of this work is two-fold. First, we show that the solution of the LD problem provides a lower bound that is as tight as that provided by the SDP relaxation.
Because this LD bound is still challenging to solve, we propose a relaxation of that LD bound, which is less tight than SDP relaxation but as tight as SOC relaxation.
Second, we derive a PBM to solve for the LD problem.
Overall, we believe that our work provides useful theoretical and computational results that enable the analysis and solution of large ACOPF problems. 

This paper is structured as follows. Section \ref{sec:net_decomp} derives the LD problem and associated network decomposition approach. Section \ref{sec:solution_methods} proposes the PBM to solve the LD problem and establishes convergence properties. Section \ref{sec:numerical} reports numerical results on standard test cases. Section \ref{sec:conclusion} concludes the paper.

\section{Lagrangian Dual Bound of ACOPF}
\label{sec:net_decomp}
We consider a rectangular ACOPF problem of the form
\begin{subequations}
\label{ACOPF}
\begin{align}
    z = \min_x \quad 
    & \textstyle\sum_{i\in\calN} f_i(p_i^G) \label{acopf:obj} \\
    \text{s.t.} \quad
    & x:=(W,p^f,p^t,q^f,q^t,p^G,q^G) \in \calF(\calG), \label{acopf:constrs} \\
    & W \in \calW := \left\{v v^T : v \in \setC^{|\calN|} \right\} \subseteq \setC^{|\calN|\times |\calN|}, \label{cons:W}
\end{align}
\end{subequations}
where $W$ is used to capture the voltage bus matrix, which is complex and can be written as $W = W^R + \ibf W^I$, where $W^R := \Re(W)$ and $W^I := \Im(W)$.
Let $\calF(\calG)$ be the feasible set defined by linear and convex quadratic constraints with respect to graph $\calG=(\calN,\calL)$, whereas the nonconvex constraints are isolated in set $\calW$.
The detailed formulation is available and adopted from~\cite{powersys}.
We assume, without loss of generality, that $f_i$ is linear for all $i \in \calN$ (higher order polynomial terms can be captured via lifting). We observe that problem \eqref{ACOPF} is equipped with linear or convex quadratic constraints except for the rank-1 constraints \eqref{cons:W}, which is nonconvex quadratic.

The proposed Lagrangian dual bound is computed by performing network decomposition and  then applying a Lagrangian relaxation of the coupling constraints between partitions. Consider a set of partitions $\{\calG_k\}_{k\in\calK}$ of $\calG$, where $\calK := \{1,\dots,K\}$ is the index set of partitions. We define each partition $\calG_k = (\calN_k,\calL_k)$ such that 
\begin{itemize}
    \item $\cup_{k\in\calK} \calN_k = \calN$ and $\calN_k \cap \calN_{k'} = \emptyset$ for $k\neq k'$ and
    \item $\calL_k := \cup_{i\in\calN_k} \left(\calL_i^f\cup\calL_i^t\right)$ and $\cup_{k\in\calK} \calL_k = \calL$,
\end{itemize}
where $\calL_i^f$ is the set of lines with the origin of bus $i$ and $\calL_i^t$ is the set of lines with the destination of bus $i$. 
Based on the partition, we also define a set of cuts that consists of lines that are shared by two graphs, by $\calC_k := \cup_{k'\in\calK\setminus\{k\}} (\calL_k \cap \calL_{k'})$. 
The set of all lines in cuts is then defined as $\calC := \cup_{k\in\calK}\calC_k$. 
Note that by definition, $\calC_k \subseteq \calL_k$. 
For each subgraph $\calG_k$, we define local decision variables $W_k$ for the voltage bus matrix and $p^f_{lk}, p^t_{lk}, q^f_{lk}, q^t_{lk}$ for each $l \in \calL_k$. 
In addition, we define a set of global decision variables $(p^f_l, p^t_l, q^f_l, q^t_l)$ for any $l\in\calC$. 

The problem~\eqref{ACOPF} can then be reformulated as
\begin{subequations}
\label{net_acopf}
\begin{align}
    \min \quad 
    & \textstyle\sum_{k\in\calK} \sum_{i\in\calN_k} f_i(p_i^G) \label{net_acopf:obj} \\
    \text{s.t.} \quad
    & (W_k,p_k^f,p_k^t,q_k^f,q_k^t,p_k^G,q_k^G) \in \calF(\calG_k), \ \forall k\in\calK, \label{net_acopf:constrs} \\
    & W_k \in \calW, \quad \forall k\in\calK, \label{net_acopf:w} \\
    & p_l^f = p_{lk}^f, \
    p_l^t = p_{lk}^t, \ \forall l\in \calC_k, \, k\in \calK, \label{net_acopf:pl} \\
    & q_l^f = q_{lk}^f, \ 
    q_l^t = q_{lk}^t, \ \forall l\in \calC_k, \, k\in \calK, \label{net_acopf:ql} \\
    & W_k = W, \quad \forall k\in\calK.
\end{align}
\end{subequations}
The objective function \eqref{net_acopf:obj} and constraints \eqref{net_acopf:constrs} and \eqref{net_acopf:w} capture the standard ACOPF formulation for each partition. Constraints \eqref{net_acopf:pl}--\eqref{net_acopf:ql} ensure that flow variables in different partitions maintain the same value for the same line in any cut. Thus, the global variables $x_k:=(W_k, p^f_l, p^t_l, q^f_l, q^t_l,p_k^G,q_k^G)$ are not specific to any partition. 

For the following discussion, we consider a simplified form of problem \eqref{net_acopf} as
\begin{subequations}
\label{simplified}
\begin{align}
    \min_{W_k,x_k,y,y_k} \quad 
    & \textstyle\sum_{k\in\calK} f_k(x_k) \\
    \text{s.t.} \quad
    & (W_k,x_k,y_k) \in \calR_k, \quad \forall k\in\calK, \\
    & y - y_k = 0, \quad (\lambda_k) \quad \forall k\in\calK, \label{simplified:linking}
\end{align}
\end{subequations}
where $x_k$ collects all local variables of partition $k$, $y_k$ collects all flow variables of partition $k$ associated with lines in cut $\calC_k$, $\calR_k := \{ (W_k,x_k,y_k) \in \calF(\calG_k): W_k \in \calW \}$ for each partition $k \in \calK$, and $y$ collects all global variables. Similarly, $f_k(x_k) := \sum_{i\in\calN_k} f_i(p_i^G)$, and $\lambda_k$ is the Lagrangian multiplier for the corresponding constraint. Note that many elements of matrix $W_k$ are free variables for a given $\calG_k$.

To decompose the problem, we relax the linking constraints \eqref{simplified:linking} and formulate the Lagrangian dual function $D(\lambda)$ as 
\begin{subequations}
\label{LDfunction}
\begin{align}
    \min_{W_k,x_k,y,y_k} \quad 
    & \textstyle\sum_{k\in\calK} \left[ f_k(x_k) - \lambda_k^T (y_k - y) \right]\\
    \text{s.t.} \quad
    & (W_k,x_k,y_k) \in \calR_k, \quad \forall k\in\calK.
\end{align}
\end{subequations}
By weak duality, $z \geq D(\lambda)$ for any $\lambda$.
To make the Lagrangian dual function bounded, we require $y^T \sum_{k\in\calK} \lambda_k = 0$.
The Lagrangian dual bound is obtained by solving the following problem:
\begin{align}
\label{LD}
    z_{LD} = \max_{\lambda_k} \left\{ \textstyle\sum_{k\in\calK} D_k(\lambda_k) : \textstyle\sum_{k\in\calK} \lambda_k = 0 \right\},
\end{align}
where the Lagrangian function $D_k(\lambda_k)$ is concave in $\lambda_k$ for all $k\in\calK$ and defined as follows:
\begin{subequations}
\label{Lagrange_subproblem}
\begin{align}
    D_k(\lambda_k) := \min_{W_k,x_k,y_k} \quad 
    & f_k(x_k) - \lambda_k^T y_k\\
    \text{s.t.} \quad
    & (W_k,x_k,y_k) \in \calR_k,
\end{align}
\end{subequations}
Note that the Lagrangian dual problem~\eqref{LD} is a concave maximization problem that may be solved by convex optimization algorithms.
Moreover, the Lagrangian subproblem \eqref{Lagrange_subproblem} is equivalent to the following convex relaxation because of the linear objective function:
\begin{subequations}
\label{Lagrange_subproblem_conv}
\begin{align}
    D_k(\lambda_k) = \min_{W_k,x_k,y_k} & \; f_k(x_k) - \lambda_k^T y_k\\
    \text{s.t.} & \; (W_k,x_k,y_k) \in \conv\left( \calR_k \right).
\end{align}
\end{subequations}
where $\conv(\cdot)$ is the convex hull opeartor: the smallest convex set that encloses set $\cdot$.
It is well known that the Lagrangian dual $z_{LD}$~\eqref{LD} equals the optimal value of the following problem (e.g., see Proposition 1 of~\cite{kim2018algorithmic}):
\begin{subequations}
\label{LDtight}
\begin{align}
    \min_{W_k,x_k,y,y_k} \
    & \textstyle\sum_{k\in\calK} f_k(x_k) \\
    \text{s.t.} \
    & (W_k,x_k,y_k) \in \conv\left( \calR_k \right) \ \forall k\in\calK, \\
    & y - y_k = 0 \ \forall k\in\calK.
\end{align}
\end{subequations}

Note that $D_k$ is still challenging to compute because each $D_k$ keeps a full copy of voltage bus matrix $W_k$ and it is still subject to rank-1 constraint. To enable more scalable computation, we reduce the size of the matrix $W_k$ to $\hat{W}_k$ that only captures the buses within subnetwork $k$, plus the neighboring buses. Furthermore, we convexify the subproblem by restricting $\hat{W}_k$ to be positive semi-definite, instead of rank-1. This leads to the following convex subproblem:
\begin{subequations}
\label{Lagrange_subproblem_sdp}
\begin{align}
    \hat{D}_k(\lambda_k) := \min_{\hat{W}_k,x_k,y_k} & \; f_k(x_k) - \lambda_k^T y_k\\
    \text{s.t.} & \; (\hat{W}_k,x_k,y_k) \in \calF(\calG_k), \;
    \hat{W}_k \succeq 0
\end{align}
\end{subequations}
Now we can apply off-the-shelf SDP solvers to each $\hat{D}_k$, which gives us a relaxed Lagrangian dual bound:
\begin{align}
\label{LD_SDP}
    \hat{z}_{SDP} := \max_{\lambda_k} \left\{ \textstyle\sum_{k\in\calK} \hat{D}_k(\lambda_k) : \sum_{k\in\calK} \lambda_k = 0 \right\}
\end{align}

We now characterize the tightness of $z_{LD}$ and $\hat{z}_{SDP}$. First, we present the SDP relaxation of~\eqref{simplified}; note that the set of positive semidefinite matrices is obtained by taking the convex hull $\conv(\calW)$.
Therefore, the SDP relaxation is given by
\begin{align}
\label{sdp}
    z^{SDP} := \min_{W_k,x_k,y,y_k} \ 
    & \textstyle\sum_{k\in\calK} f_k(x_k) \\
    \text{s.t.} \
    & (W_k,x_k,y_k) \in \calF(\calG_k), \ \forall k\in\calK, \notag \\
    & W_k \in \conv(\calW), \ \forall k\in\calK, \notag \\
    & y - y_k = 0, \ \forall k\in\calK. \notag 
\end{align}
Note that \eqref{sdp} is equivalent to the SDP relaxation of~\eqref{ACOPF}.

Now we present the main theoretical result of this work: the network decomposition attains a tighter lower bound than SDP relaxation.
\begin{theorem} 
$z_{LD} \geq z^{SDP} \geq \hat{z}_{SDP} \geq z^{SOC}$.
\end{theorem}
\begin{IEEEproof}
The inequalities $z^{SDP} \geq \hat{z}_{SDP} \geq z^{SOC}$ hold because a positive semi-definite matrix implies positive semi-definiteness for all submatrices, and positive semi-definite submatrices ensure that the SOC constraints are always satisfied.

Now we establish $z_{LD} \geq z^{SDP}$. Denote by $\calD$ and $\calD^{SDP}$, respectively, the feasible region of \eqref{LDtight} and \eqref{sdp}. For any $(W_k, x_k, y, y_k) \in \calD$, it suffices to show $(W_k, x_k, y, y_k) \in \calD^{SDP}$. First, $y = y_k$ is satisfied for all $k \in \calK$ because of the feasibility in \eqref{LDtight}. Since $(W_k, x_k, y_k) \in \conv\left( \calR_k \right)$ for all $k\in\calK$, there exist, for any $k \in \calK$, finitely many $(W_k^1, x_k^1, y^1, y_k^1), ..., (W_k^n, x_k^n, y^n, y_k^n)$, and scalars $\alpha^k_1, ..., \alpha^k_n$ such that
\begin{subequations}
\label{convex_comb}
\begin{gather}
    \textstyle\sum_{i = 1}^n \alpha^k_i = 1, \ \alpha^k_i \geq 0 \quad \forall i = 1,...,n \label{convex_comb_1} \\
    (W_k, x_k, y_k) = \textstyle\sum_{i = 1}^n \alpha^k_i \left( W_k^i, x_k^i, y_k^i \right) \label{convex_comb_4} \\
    (W_k^i, x_k^i, y_k^i) \in \calF(\calG_k), \ W_k^i \in \calW, \ \forall i = 1,...,n. \label{convex_comb_5}
\end{gather}
\end{subequations}
Equations \eqref{convex_comb_1}--\eqref{convex_comb_4} imply that $W_k \in \conv(\calW)$ for all $k \in \calK$. Furthermore, from \eqref{convex_comb}, we have $(W_k, x_k, y_k) \in \calF(\calG_k)$ for any $k \in \calK$ because $\calF(\calG_k)$ is convex, since it is defined by linear and convex quadratic constraints.
\end{IEEEproof}

\begin{figure}[!htb]
    \centering
    \includegraphics[width=0.35\textwidth]{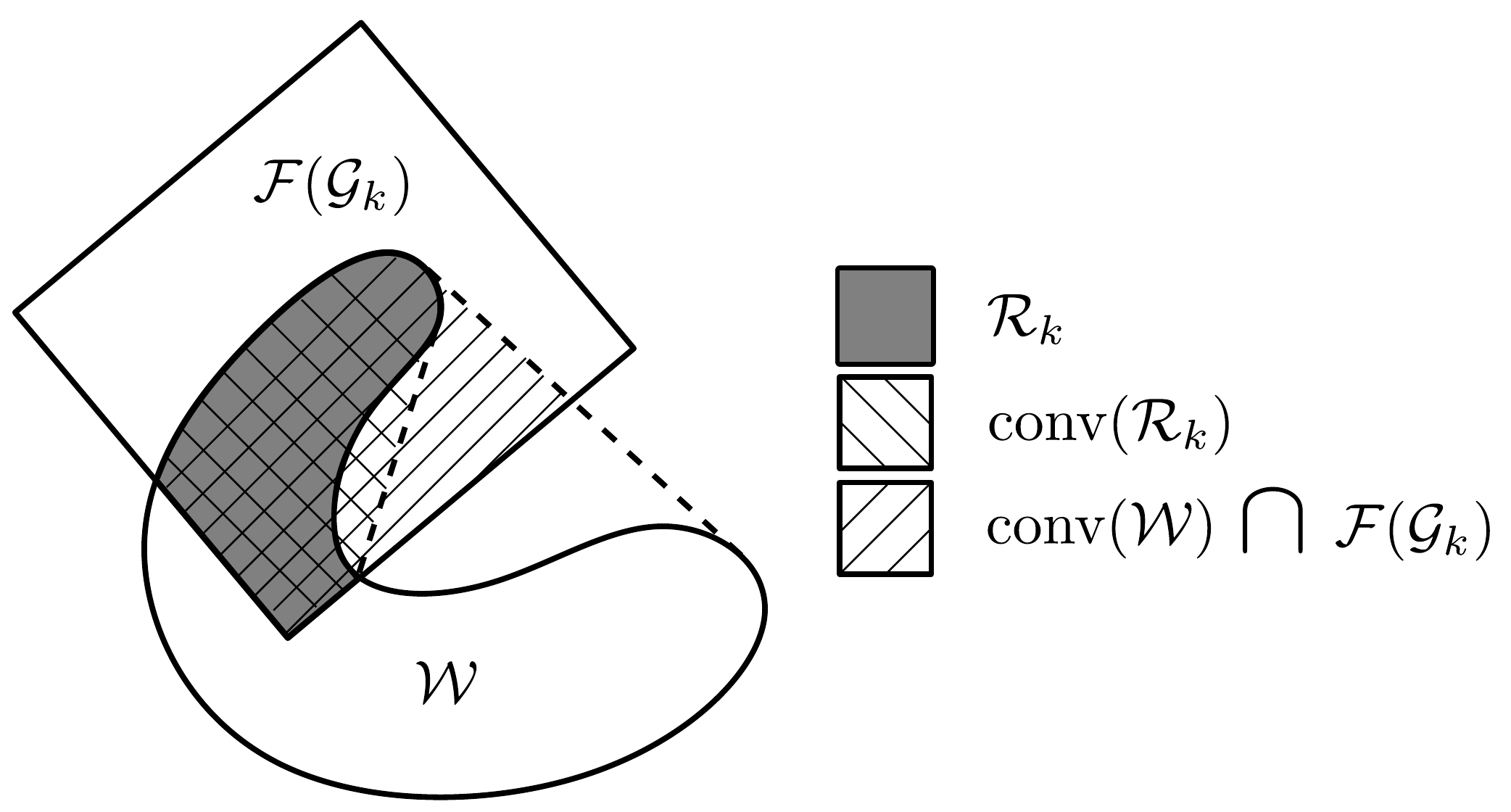}
    \caption{Comparison between feasible sets of SDP relaxation and network decomposition.}
    \label{fig:tight_bound}
\end{figure}

Figure \ref{fig:tight_bound} gives a graphical illustration for the first inequality $z_{LD} \geq z^{SDP}$. Note that this argument relies on the assumption of linear objective function, which enables the equivalence between \eqref{Lagrange_subproblem} and \eqref{Lagrange_subproblem_conv}. Furthermore, if at the optimal dual $\lambda$, the solution for the subproblems $(W_k, x_k, y_k)$ satisfy all the coupling constraints, then the point $(W_k, x_k, y, y_k)$ is feasible for the primal ACOPF and is globally optimal ($z_{LD} = z$). 

\section{Computation of Lagrangian Dual Bound}
\label{sec:solution_methods}

Now we show how to solve \eqref{LD_SDP} to obtain $\hat{z}_{SDP}$. This can be achieved by using a proximal bundle method (PBM), following the algorithmic steps in \cite{kiwiel1990proximity}. Here we describe the PBM algorithm using notation from our network decomposition framework. The method formulates a model function to search for the optimal $\lambda$ values by taking the linear outer approximations of the Lagrangian dual functions $D_k$. 

Let $r, l \in \mathbb{Z}_+$ be the indices for major and minor iterations, respectively. The model function is formulated as follows:
\begin{subequations}
\begin{align}
    m_{r,l}(\lambda) := \max_{\theta_k} & \quad \textstyle\sum_{k \in \calK} \theta_k \\
    \text{s.t.} & \quad \theta_k \leq \hat{D}_k(\lambda^i) - (y^i_k)^T(\lambda-\lambda^i), \nonumber \\
    & \qquad\qquad i \in \calB^{r,l}, k \in \calK,
\end{align}
\end{subequations}
where $\calB^{r,l}$ is the set of cut indices at iteration $(r,l)$. We note that $-y^i_k \in \partial \hat{D}_k(\lambda^i)$, and $m_{r,l}$ outer-approximates $\hat{D}$; that is, $m_{r,l}(\lambda) \geq \hat{D}(\lambda)$  for any $\lambda$. 

At each iteration $(r,l)$ we solve the master problem
\begin{equation}
\label{proximal_master}
    \max_\lambda \ m_{r,l}(\lambda) - u_{r,l}\|\lambda - \lambda^r\|^2 / 2,
\end{equation}
where $\lambda^r$ is the starting point for all minor iterations under major iteration $r$ and  $u_{r,l} > 0$ is the proximal weight.
At each minor iteration $l$, the master problem finds a new trial point $\lambda^{r,l}$ for evaluating $\hat{D}(\cdot)$. The predicted increase of $\hat{D}(\cdot)$ at the trial point is defined as
\begin{equation}
    v_{r,l} := m_{r,l}(\lambda^{r,l}) - \hat{D}(\lambda^r).
\end{equation}
The sequence of minor iterations under one major iteration is terminated if the sufficient decrease condition is satisfied:
\begin{equation}
    \label{eq:suff_decrease}
    \hat{D}(\lambda^{r,l}) \geq \hat{D}(\lambda^r) + m_L v_{r,l},
\end{equation}
where $m_L \in (0, 1/2)$. This is referred to as a serious step. At a serious step, the algorithm updates the starting center $\lambda^{r+1} \leftarrow \lambda^{r,l}$ and enters the next major iteration. Otherwise, if the sufficient decrease condition is not satisfied, the algorithm executes a null step, in which the proximal weight is possibly increased to restrict search range. The algorithm terminates if $v_{r,l} \leq \epsilon \cdot (1 + |\hat{D}(\lambda^r)|)$
for some $\epsilon \in \mathbb{R}_+$. We summarize the algorithmic steps in Algorithm~\ref{prox_algo}. 

\begin{algorithm}[h!]
\SetAlgoLined
 Initialize $\lambda^0$, $u > 0$, $m_L \in (0, 1/2)$, $\epsilon > 0$, $\calB^{0,0} \leftarrow \{0\}$, $r \leftarrow 0$, $l \leftarrow 0$. \\
 Solve \eqref{Lagrange_subproblem_sdp} to find $\hat{D}_k(\lambda^0)$ and $y_k^0$ for $k \in \calK$. \\
 Initialize the model function $m_{0,0}$. \\
 $v_{r,l} \leftarrow m_{r,l}(\lambda^{r,l}) - \hat{D}(\lambda^r)$. \\
 \While{$v_{r,l} \leq \epsilon (1 + |\hat{D}(\lambda^r)|)$}{
  Solve the proximal master \eqref{proximal_master} for $\lambda^{r,l}$. \\
  Solve the Lagrangian dual \eqref{Lagrange_subproblem_sdp} for $\hat{D}_k(\lambda^{r,l})$ and $y^{r,l}_k$ for $k \in \calK$. \\
  $v_{r,l} \leftarrow m_{r,l}(\lambda^{r,l}) - \hat{D}(\lambda^r)$. \\
  Add to $\calB^{r,l}$ cuts for all $k \in \calK$ $\theta_k \leq \hat{D}_k(\lambda^{r,l}) - (y^{r,l}_k)^T(\lambda-\lambda^{r,l})$ \\
  \uIf{$\hat{D}(\lambda^{r,l}) - \hat{D}(\lambda^r) \geq m_Lv_{r,l}$}{
        $\lambda^{r+1} \leftarrow \lambda^{r,l}, \calB^{r+1, 0} \leftarrow \calB^{r,l}$ \\
        Update proximal weight $u^{r+1,0}$ \\
        $r \leftarrow r+1$, $l \leftarrow 0$ \\
   }
  \Else{
        $\calB^{r,l+1} \leftarrow \calB^{r,l}$, and
        choose $u^{r,l+1} \in [u^{r,l}, \infty)$ \\
        $i \leftarrow \min\{i-1, -1\}$, $l \leftarrow l + 1$ \\
  }
 }
 \caption{Proximal Bundle Method}\label{prox_algo}
\end{algorithm}

The work of \cite{kim2019asynchronous} has established that the trust-region variant of the proximal bundle method converges to the dual optimality and terminates with a finite number of iterations, provided that the Lagrangian subproblems $\hat{D}_k$ are solved to global optimality. 

\section{Numerical Experiments}
\label{sec:numerical}

In this section we demonstrate the numerical performance of the dual decomposition of ACOPF.
We have implemented parallel network decomposition algorithms in {\tt NetDecOPF.jl}, in which subnetwork ACOPF formulations are generated by using {\tt PowerModels.jl}~\cite{powermodels}. 
The decomposition of the network is performed by using the METIS library through the {\tt Metis.jl} interface, and the subproblems are solved by {\tt MOSEK v9.2.24}.
To avoid the slow progress status of {\tt Mosek}, we relax the convergence parameters for primal and dual feasibilities and relative gap to $10^{-4}$.
The master problem is solved by OSQP~\cite{osqp} compiled with Intel MKL Pardiso for linear solver.
We initialize the Lagrangian multiplier $\lambda^0$ to zeros and the termination tolerance $\epsilon$ to $10^{-4}$ in Algorithm~\ref{prox_algo}.

We use the test instances whose SOC gaps are known to be larger than 5\% among the instances from the typical operating conditions and congested operating conditions in PGLib-OPF (v20.07). 
All experiments were conducted using Julia and run on the Bebop cluster at Argonne National Laboratory. 
To avoid any racing conditions, each run used up to 32 cores only per compute node to parallelize each subnetwork problem solution by the network decomposition. 
The subnetwork problems were distributed in a round-robin fashion such that each problem is assigned to one core.

\begin{table*}[ht!]
    \caption{Computational performance of the Lagrangian dual of ACOPF for different test problems (using MOSEK).}
    \label{table:results}
    \begin{center}
    \begin{tabular}{l|r|r|r|rrrr}
     \hline
          &       &       &       &  \multicolumn{4}{c}{Gap (\%)} \\
     Case & Nodes & Edges & $|\calN|$       & $z_{QC}$ & $z_{SOC}$ & $z_{SDP}$ & $\hat{z}_{SDP}$ \\
     \hline
     case5\_pjm                & 5     & 6     & 2  & 14.55 & 14.55 & 5.22 & 5.31 \\
     case14\_ieee\_\_api       & 14    & 20    & 2  & 5.13  & 5.13  & 0.02 & 1.14 \\
     case24\_ieee\_rts\_\_api  & 24    & 38    & 3  & 13.01 & 17.88 & 2.07 & 14.02 \\
     case30\_as\_\_api         & 30    & 41    & 3 & 44.61 & 44.61 & * & 38.86 \\
     case30\_ieee              & 30    & 41    & 3 & 18.81 & 18.84 & 0.02 & 1.45 \\
     case30\_ieee\_\_api       & 30    & 41    & 3 & 5.46  & 5.46  & * & 1.85 \\
     case73\_ieee\_rts\_\_api  & 73    & 120   & 7 & 11.05 & 12.87 & 2.91 & 4.64 \\
     case89\_pegase\_\_api     & 89    & 210   & 8 & 23.07 & 23.11 & 21.91 & 23.56 \\
     case118\_ieee\_\_api      & 118   & 186   & 11 & 29.74 & 29.97 & 11.81 & 26.14\\
     case179\_goc\_\_api       & 179   & 263   & 17 & 5.93 & 9.88 & 0.56 & 0.67 \\
     \hline
    \end{tabular}
    \end{center}
    \footnotesize{* denotes the instances that are solved to \texttt{UNKNOWN\_RESULT\_STATUS} for both primal and dual status. }
    \vspace{-0.1in}
\end{table*}

Table~\ref{table:results} reports the duality gaps obtained by the network decomposition, as well as the gaps of QC, SOC and SDP relaxations. 
Each system is partitioned by different number of partitions to ensure each subproblem has around 10 buses for large systems. 
The QC and SOC bounds were obtained from the Github repository of PGLib-OPF (v20.07).
The sparsity-exploiting SDP relaxation was formulated by using \texttt{SparseSDPWRMPowerModel} in \texttt{PowerModels.jl} and solved by using \texttt{MOSEK v9.2.24}. 
The gaps were calculated with the AC objective values reported in PGLib-OPF. 


We observe that our network decomposition obtains a bound between SDP and SOC, consistent with our theoretical result. 
The only exception comes from {\tt case89\_pegase\_\_api}, where the gap between SOC and SDP is relatively small. 
This is possible in practice because we are using relaxed convergence parameters for {\tt Mosek} when solving for SDP subproblems (so the Lagrangian dual bound obtained is lower than the actual value). 
Despite the relaxation, we are still able to achieve theoretical bound order for most cases.
We also highlight that the SDP relaxation models for two systems could not be solved to good result status due to numerical issues.

\begin{figure}[htb!]
    \centering
    \includegraphics[width=0.4\textwidth]{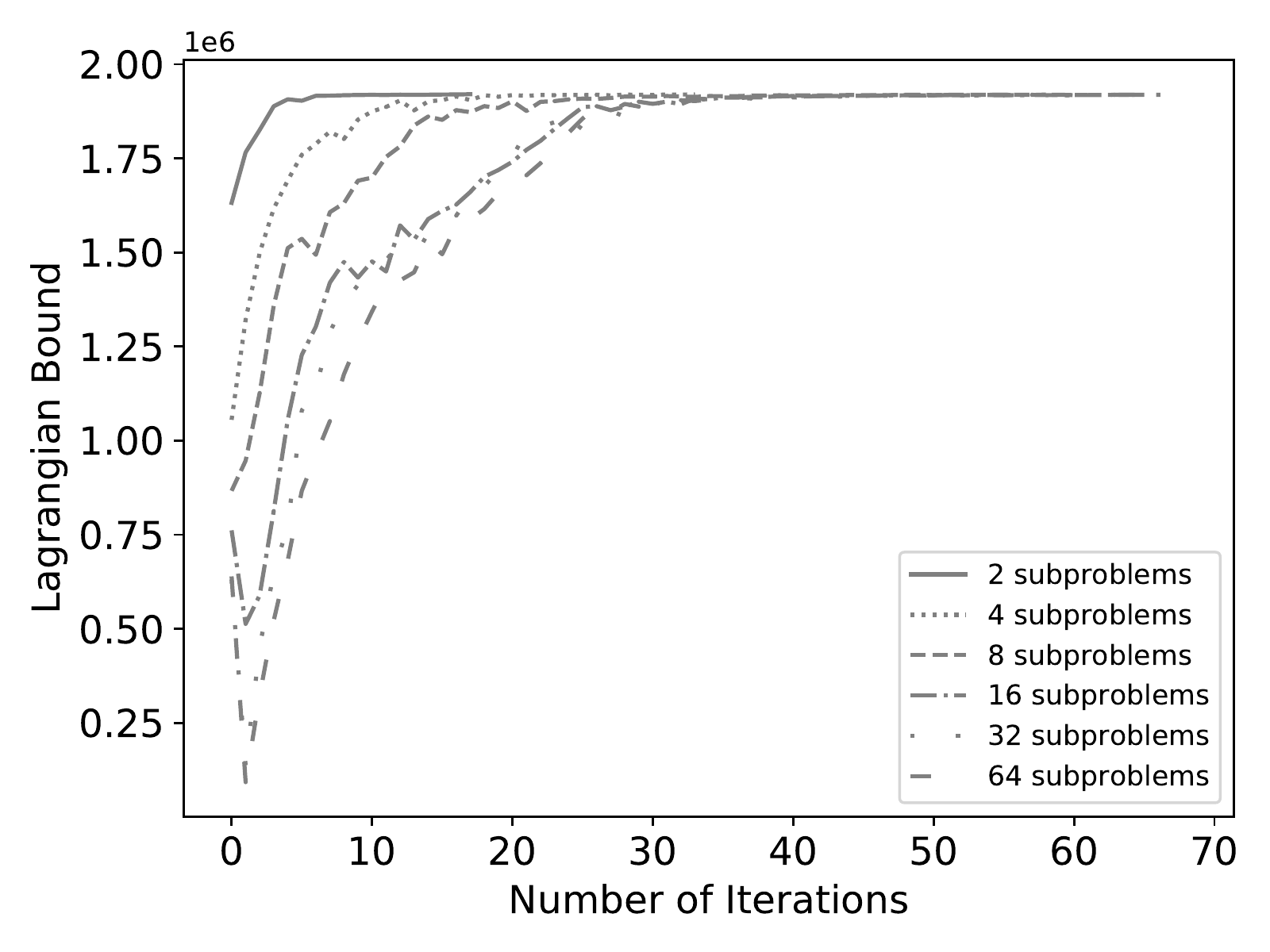}
    \caption{Iteration of the Lagrangian dual decomposition of {\tt case179\_goc\_\_api} with different number of subproblems}
    \label{fig:scale_iter}
\end{figure}

The Lagrangian dual method provides a valid lower bound for any given dual multiplier at any iteration.
This property can be useful in practice when computational resources (e.g., time, CPUs) are limited. 
In particular, we observe that the Lagrangian dual bound is already very tight at the very first iteration when the problem is decomposed into a smaller number of subproblems (see Figure~\ref{fig:scale_iter}).

The decomposition to a large number of subproblems decreases the Lagrangian dual bound, as shown in Figure~\ref{fig:scale_iter}.
However, we emphasize that even with many subproblems, the Lagrangian dual bound should still be as tight as that of the SOC relaxation.
In addition, we note that the Lagrangian dual decomposition provides a lower bound for all iterations; the method is simply searching for the optimal dual values that provide the tightest bound. 
For practical purposes, one can run the first iteration only and get a valid lower bound. However, there will be no guarantee for the quality of the bound. 

\begin{figure}[htb!]
    \centering
    \includegraphics[width=0.4\textwidth]{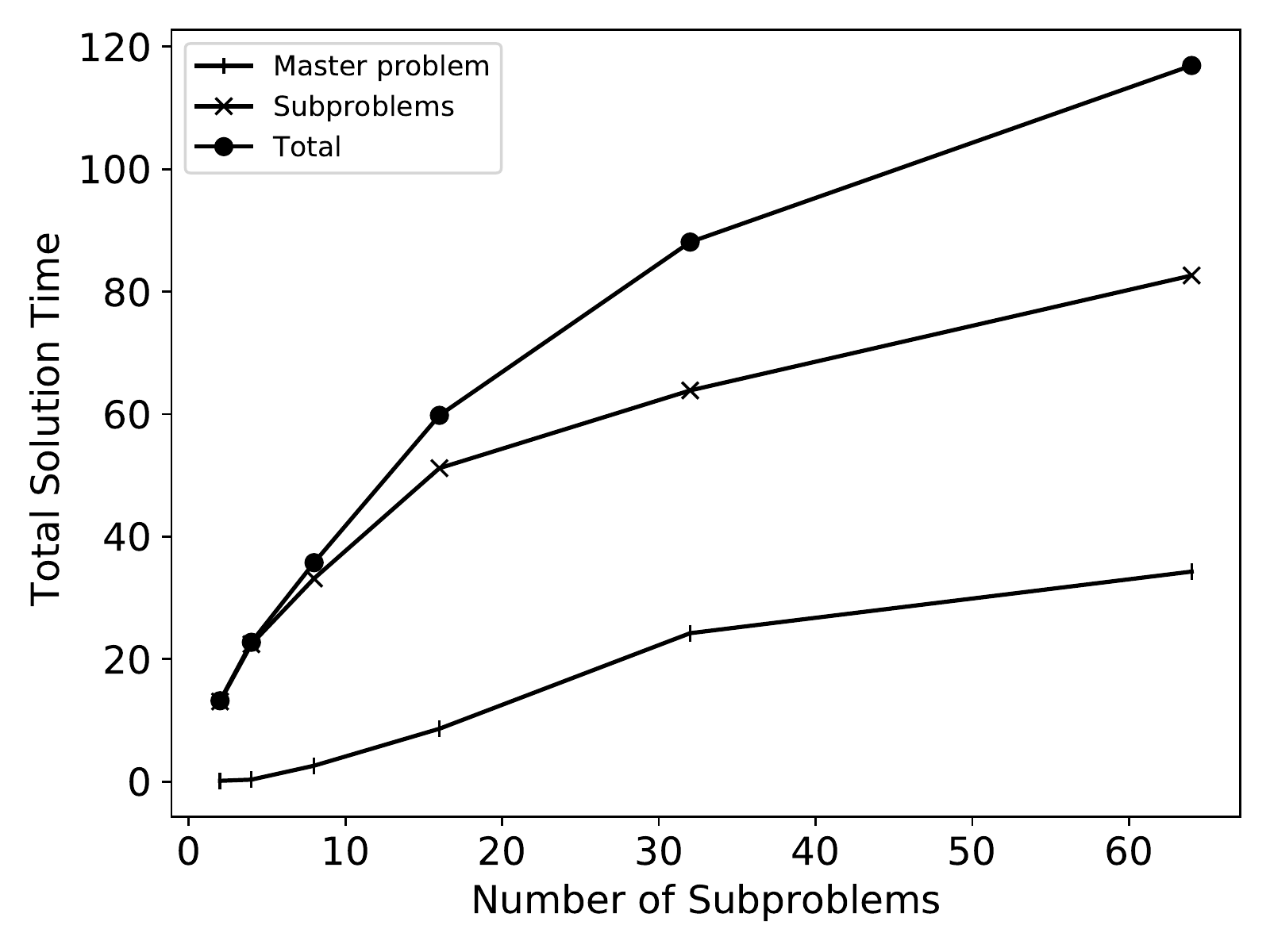}
    \caption{Solution times of the Lagrangian dual decomposition of {\tt case179\_goc\_\_api} with different number of subproblems}
    \label{fig:scale_time}
    \vspace{-0.1in}
\end{figure}

We conclude this section by discussing a major computational challenge observed in our numerical study.
Figure~\ref{fig:scale_time} shows that the changes of the master problem solution time, parallel subproblem solution time, and the total solution time as the number of subproblems increases.
We observe that unlike the scenario decomposition of stochastic programming (e.g.,~\cite{kim2018algorithmic}), the master problem of the network decomposition is extremely difficult to solve and time-consuming, particularly with a large number of subproblems.
The reasons may be that the master problem of the network decomposition is less structured than that of the stochastic programming and also that each iteration requires to add multiple dense inequality constraints to the master problem, as also discussed in temporal decomposition~\cite{kim2018temporal}.
We also observe that the solution time of master problem increases with number of partitions, but is still dominated by the subproblem solution time. 
However, for larger cases not shown here, we observe that master problem starts to become a serious bottleneck.

\section{Concluding Remarks and Future Work}
\label{sec:conclusion}
We propose a Lagrangian dual bound for large-scale ACOPF problems and prove its theoretical tightness compared with SDP and SOC. We show that the bound can be computed by applying a PBM; notably, this approach can employ parallel computation of the subproblems. 
Our numerical experiments demonstrate the theoretical tightness properties.
We also identified and discussed the computational bottleneck of the proposed PBM.
We will seek for a scalable method for solving the master problem as future work. 

\bibliographystyle{IEEEtran}
\bibliography{refs}

\end{document}